\documentclass[10pt,draft]{amsart}
\usepackage{amsmath,amssymb,amsthm}
\begin{document}
\newtheorem{lem}{Lemma}[section]
\newtheorem{prop}{Proposition}[section]
\newtheorem{cor}{Corollary}[section]
\numberwithin{equation}{section}
\newtheorem{thm}{Theorem}[section]

\theoremstyle{remark}
\newtheorem{example}{Example}[section]
\newtheorem*{ack}{Acknowledgments}

\theoremstyle{definition}
\newtheorem{definition}{Definition}[section]

\theoremstyle{remark}
\newtheorem*{notation}{Notation}
\theoremstyle{remark}
\newtheorem{remark}{Remark}[section]

\newenvironment{Abstract}
{\begin{center}\textbf{\footnotesize{Abstract}}%
\end{center} \begin{quote}\begin{footnotesize}}
{\end{footnotesize}\end{quote}\bigskip}
\newenvironment{nome}

{\begin{center}\textbf{{}}%
\end{center} \begin{quote}\end{quote}\bigskip}

\newcommand{\triple}[1]{{|\!|\!|#1|\!|\!|}}

\newcommand{\xx}{\langle x\rangle}
\newcommand{\ep}{\varepsilon}
\newcommand{\al}{\alpha}
\newcommand{\be}{\beta}
\newcommand{\de}{\partial}
\newcommand{\la}{\lambda}
\newcommand{\La}{\Lambda}
\newcommand{\ga}{\gamma}
\newcommand{\del}{\delta}
\newcommand{\Del}{\Delta}
\newcommand{\sig}{\sigma}
\newcommand{\ome}{\omega}
\newcommand{\Ome}{\Omega}
\newcommand{\C}{{\mathbb C}}
\newcommand{\N}{{\mathbb N}}
\newcommand{\Z}{{\mathbb Z}}
\newcommand{\R}{{\mathbb R}}
\newcommand{\Rn}{{\mathbb R}^{n}}
\newcommand{\Rnu}{{\mathbb R}^{n+1}_{+}}
\newcommand{\Cn}{{\mathbb C}^{n}}
\newcommand{\spt}{\,\mathrm{supp}\,}
\newcommand{\Lin}{\mathcal{L}}
\newcommand{\SSS}{\mathcal{S}}
\newcommand{\F}{\mathcal{F}}
\newcommand{\xxi}{\langle\xi\rangle}
\newcommand{\eei}{\langle\eta\rangle}
\newcommand{\xei}{\langle\xi-\eta\rangle}
\newcommand{\yy}{\langle y\rangle}
\newcommand{\dint}{\int\!\!\int}
\newcommand{\hatp}{\widehat\psi}
\renewcommand{\Re}{\;\mathrm{Re}\;}
\renewcommand{\Im}{\;\mathrm{Im}\;}

\title
[Radiality of minimizers to S-P-S energy]{On the radiality of constrained minimizers to the
Schr\"odinger-Poisson-Slater energy}

\author{Vladimir Georgiev}
\address{Vladimir Georgiev: Universit\`a di Pisa,
Dipartimento di matematica, Largo B. Pontecorvo 5, 56124, Pisa, Italy}
\email{georgiev@dm.unipi.it}

\author{Francesca Prinari}
\address{Francesca Prinari: Universit\`a di Ferrara,
Dipartimento di Matematica, via Machiavelli 35, 44100, Ferrara, Italy}
\email{francesca.prinari@unife.it}

\author{Nicola Visciglia}
\address{Nicola Visciglia: Universit\`a di Pisa, Dipartimento di Matematica, Largo B.
Pontecorvo 5, 56100 Pisa, Italy}
\email{viscigli@dm.unipi.it}

%\subjclass[2000]{35J10, 35L05.}

%\keywords{Fourier restriction theorems, Strichartz estimates}

\maketitle

\maketitle

\begin{abstract}
We study the radial symmetry
of minimizers to the Schr\"odinger-Poisson-Slater
(S-P-S) energy:
$$\inf_{\substack{u\in H^1(\R^3)\\
\|u\|_{L^2(\R^3)}=\rho}}
\frac 12 \int_{\R^3} |\nabla u|^2 + \frac 14
\int_{\R^3} \int_{\R^3} \frac{|u(x)|^2|u(y)|^2}{|x-y|}dxdy
- \frac 1p \int_{\R^3} |u|^p dx$$
provided that  $2<p<3$ and $\rho$ is small. The main result shows that minimizers are radially symmetric modulo suitable translation.
\end{abstract}

The following minimization problem associated to Schr\"odinger-Poisson-Slater
(S-P-S)
energy functional has been extensively studied in the literature
(see  for instance \cite{BS1}, \cite{BS}, \cite{CattoL}, \cite{SS} and
all the references therein):
\begin{equation}\label{I}
I_{\rho, p}=\inf_{\substack{u\in H^1(\R^3)\\
\|u\|_{L^2(\R^3)}=\rho}} {\mathcal E}_{p}(u)
\end{equation}
where
$${\mathcal E}_{p}(u)=
\frac 12 \int_{\R^3} |\nabla u|^2 + \frac 14
\int_{\R^3} \int_{\R^3} \frac{|u(x)|^2|u(y)|^2}{|x-y|}dxdy
- \frac 1p \int_{\R^3} |u|^p dx.$$
The corresponding set of minimizers will be denoted since now on by $\mathcal M_{\rho, p}$.
It has been proved in \cite{SS} (based on the technique introduced in \cite{CattoL})
that ${\mathcal M}_{8/3, \rho}\neq \emptyset$ provided that
$0<\rho<\rho_0$ for a suitable $\rho_0>0$ (i.e. under a smallness assumption on the charge).
In \cite{BS} it is proved that $\mathcal M_{\rho, p}\neq \emptyset$
provided that $\rho>0$ is small and
$2<p<3$. In \cite{BS1} it is treated the case $3<p< \frac{10}3$
and $\rho$ sufficiently large.\\
The main aim to look at the minimization problem
\eqref{I} is to construct (following the original argument by \cite{CL})
orbitally stable standing wave solutions to the following evolution problem
$${\bf i} \partial_t \psi + \Delta \psi- \Big (\frac 1{|x|}*|\psi|^2\Big ) \psi + \psi|\psi|^{p-2}=0
\hbox{ } (t,x)\in \R\times \R^3.$$
For the sake of completeness we recall that standing waves are solutions of the following type
$$\psi(t,x)=e^{{\bf i}\omega t} v(x)$$
for a suitable $\omega\in \R$ and $v(x)\in H^1(\R^3)$.\\\\
In this paper we study the radiality (up to translation) of the functions
in ${\mathcal M}_{\rho,p}$
provided that $\rho>0$ is small enough and $2<p<3$.

There are  different results on the symmetry of the minimizers.
The basic result due to Gidas, Ni and Nirenberg \cite{GNN} implies
the radial symmetry of the minimizers associated with the
semilinear elliptic equation
$$ \Delta u + f(u) = 0,$$ provided suitable assumptions on the function $f(u)$
are satisfied and the scalar function $u$ is positive. As in the
previous result due to Serrin \cite{S}, the proof is based on the
maximum principle and the Hopf's lemma.

The symmetry of the energy functional (even with constraint
conditions) can not imply in general the radial symmetry of the
minimizers. This phenomena was discovered and studied in the works
\cite{CZE}, \cite{E1} and \cite{E2} in the scalar case.

Different techniques have been developed in the literature
to prove the radiality
of minimizers to suitable variational problems.
We quote some of them (see also all the references therein):
\cite{BBL} where it is proved a very general radiality result
for non-negative critical points of suitable
variational problems (however Hartree type nonlinearity is not allowed),
\cite{L}, \cite{LM} where the case of
nonlocal Hartree type nonlinearity
is treated.
However, as far as we can see, those techniques
do not work in our context since the potential energy in SPS
is refocussing on the nonlocal term (the Hartree nonlinearity)
and focusing on the local term (the $L^p$ norm).\\
\\
To underline the difficulty notice that it is not obvious
to answer to the following weaker question:\\
\\
\centerline {\em  Is there at least a radially symmetric
function belonging to ${\mathcal M}_{\rho,p}$?}\\
\\
A general tool that could be useful to provide an answer to the question
above is the Schwartz
rearrangement map $u\rightarrow u^*$.
The following properties are well--known (see \cite{LL}):
$$\|\nabla u^*\|_{L^2(\R^3)}\leq \|\nabla u\|_{L^2(\R^3)}^2;$$
$$\int_{\R^3} \int_{\R^3} \frac{|u(x)|^2|u(y)|^2}{|x-y|}dxdy \geq
\int_{\R^3} \int_{\R^3} \frac{|u^*(x)|^2|u^*(y)|^2}{|x-y|}dxdy;$$
$$\|u^*\|_{L^q(\R^3)}=\|u\|_{L^q(\R^3)}.$$
As a consequence there is a competition between the kinetic energy
and the nonlocal energy which makes unclear whether or not
the set ${\mathcal M}_{\rho,p}$
is invariant under the map
$u\rightarrow u^*$ (and hence it makes useless the rearrangement
technique to provide an answer to the question raised above). \\
\\
Next we state the main result of this paper.
\begin{thm}\label{main}
For every $2<p<3$
there exists $\rho_0=\rho_0(p)>0$ such that
$$\forall (v,\rho)\in {\mathcal M}_{\rho,p}\times (0,\rho_0)
\hbox{ } \exists \tau\in \R^3
\hbox{ such that } v(x+\tau)=v(|x|+\tau) \hbox{ } \forall x\in \R^3.$$
\end{thm}
\begin{remark}
Recall that in \cite{BS} it is proved that
${\mathcal M}_{\rho,p}\neq \emptyset$
for $2<p<3$ and $\rho>0$ small.
\end{remark}
\begin{remark}
Notice that in Theorem \ref{main} the physically relevant case $p=8/3$ is allowed.
\end{remark}
Next we fix some notations.
\begin{notation}
We shall denote by $Q(x)$ the unique function such that:
\begin{enumerate} \label{qp}
 \item $Q\in H^1(\R^3)$;
\item $Q(x)$ is radially symmetric;
\item $Q(x)>0$ for every $x\in \R^3$;
\item $\|Q\|_{L^2(\R^3)}^2=1$;
\item $Q(x)$ solves the
following  elliptic problem
$$-\Delta Q + \omega_0 Q=Q|Q|^{p-2} \hbox{ on } \R^3$$
for a suitable $\omega_0>0$ (which is unique).
\end{enumerate}
(We recall that the existence and uniqueness of a function
$Q$ that satisfies the properties mentioned
above follows by combining the results in
\cite{CL},\cite{GNN}, \cite{K} provided that $2<p<\frac {10}3$).\\
If $Q(x)$ is the radial function, satisfying the above relations, then we introduce
$${\mathcal G}=\{Q(x+\tau)|\tau\in \R^3\}$$ and
for every $\tau \in \R^3$ we write $Q_\tau=Q(x+\tau)$.\\
$T_{Q}$ (resp. $T_{Q_\tau}$) denotes
the tangent space of the manifold ${\mathcal G}$
at the point $Q$ (resp. $Q_\tau$). We also denote by $T_{Q_\tau}^\perp$
the intersection of $H^1(\R^3)$ with
the orthogonal space (w.r.t. the $L^2$ scalar product)
of $T_{Q_\tau}$.\\
Let $M$ be a vector space then $\pi_M$ denotes
the orthogonal projection, with respect to the $L^2(\R^3)$ scalar product,
on the vector space $M$.\\
$H^1$ is the usual Sobolev space
endowed with the following Hilbert norm
$$\|u\|_{H^1}^2=\int_{\R^3} |\nabla u|^2 dx+\omega_0 \int_{\R^3} |u|^2 dx,$$
where $\omega_0$ is the constant introduced above.\\
$H^1_{rad}$ denotes the functions in $H^1$
that are radially symmetric.\\
$L^p$ will denote the space $L^p(\R^3)$.\\
In general $\int ... dx$ and $\int\int...dxdy$
denote $\int_{\R^3}...dx$ and $ \int_{\R^3}\int_{\R^3} ... dxdy$.\\
Assume $({\mathcal H}, (.,.))$ is an Hilbert space and 
${\mathcal F}: {\mathcal H}\rightarrow \R$
is a differentiable functional,then $\nabla_u {\mathcal F}$
is the gradient of $\mathcal F$ at the point $u\in \mathcal H$.\\
Let $(X, \|.\|)$ be a Banach space,
then $B_X(x,r)$ denotes the ball of radius $r>0$ centered in $x\in X$.\\
Let $\Phi$ be a differentiable map between two Banach
spaces $(X, \|.\|_X)$ and $(Y, \|.\|_Y)$
then $d\Phi_x\in {\mathcal L}(X,Y)$
denotes the differential of $\Phi$ at the point $x\in X$.
\end{notation}

\section{An equivalent problem}

By the rescaling
$u_\rho(x)=\rho^{\frac{4}{4-3(p-2)}} u\left
(\rho^{\frac{2(p-2)}{4-3(p-2)}}x\right )$
it is easy to check that the minimization problem
\eqref{I} is equivalent to the following one:
$$J_{\rho,p}=\inf_{\substack{u\in H^1\\
\|u\|_{L^2}=1}} \frac 12 \int |\nabla u|^2 dx +
\rho^{\alpha(p)}
\int
\int \frac{|u(x)|^2|u(y)|^2}{|x-y|}dxdy
- \frac 1p \int |u|^p dx$$
where
$\alpha(p)=\frac{16+2(p-2)-12(p-2)-4p+6(p-2)}{4-3(p-2)}$.
Notice that $\alpha(p)>0$ provided that
$2<p<3$.
Motivated by this fact we introduce the following  minimization problem
$$K_{\rho, p}=\inf_{\substack{u\in H^1\\
\|u\|_{L^2}=1}} {\mathcal E}_{\rho,p}(u)$$
where
$$ {\mathcal E}_{\rho,p}(u)=\frac 12 \int |\nabla u|^2 dx +
\rho
\int \int\frac{|u(x)|^2|u(y)|^2}{|x-y|}dxdy
- \frac 1p \int |u|^p dx.$$
We also denote by  ${\mathcal N}_{\rho, p}$ the corresponding minimizers:
$${\mathcal N}_{\rho, p}=\{v\in H^1| {\mathcal E}_{\rho,p}(v)=K_{\rho, p}\}.$$
It is easy to prove that Theorem \ref{main} is equivalent to the following proposition.
\begin{prop}\label{equiv} For every $2<p<3$
there exists $\rho_0=\rho_0(p)>0$ such that any function $v\in {\mathcal N}_{\rho, p}$
is (up to translation)
radially symmetric provided that $0<\rho<\rho_0$.
\end{prop}
The rest of the paper is devoted to the proof
of Proposition \ref{equiv}.\\\\
In the sequel the function $Q(x)$ and the constant $\omega_0>0$
are the ones defined in the introduction.\\
Next result will be useful in the sequel.
\begin{prop}\label{compactness}
Let $2<p<3$ and
$v_k\in {\mathcal N}_{\rho_k, p}$
where $\lim_{k\rightarrow \infty} \rho_k=0$.
Then up to subsequence there exists
$\tau_k\in \R^3$ such that
$$v_k(x+\tau_k)\rightarrow Q \hbox{ in } H^1.$$
\end{prop}

{\bf Proof.}
\\
{\em First step: $K_{\rho_k, p}\rightarrow K_{0, p}$
as $k\rightarrow \infty$}
\\
\\
First notice that
$$K_{0,p}\leq K_{\rho_k,p}$$
due to the positivity of $\rho_k$.
Hence it is sufficient to prove
$\limsup_{k\rightarrow \infty} K_{\rho_k,p}\leq K_{0,p}$.
This fact follows from
$$ K_{\rho_k,p}
\leq {\mathcal E}_{\rho_k,p}(Q)=
{\mathcal E}_{p}(Q)+\rho_k
\int \int \frac{|Q(x)|^2|Q(y)|^2}{|x-y|}dxdy
= K_{0,p}+ o(1).$$
\\
\\
{\em Second step:  $v_k$ converge to $Q$ up to subsequence and traslation}
\\
\\
By the previous step we deduce that
$\{v_k\}$ is a minimizing sequence for
$K_{0,p}$. As a consequence of the results proved in
\cite{CL}, \cite{GNN}, \cite{K} we deduce
that $\{v_k\}$ converge strongly
(up to translation) to $Q(x)$.

\hfill$\Box$

In next result we get a qualitative information
on the lagrange multipliers associated to the constrained minimizers
belonging to ${\mathcal N}_{\rho,p}$
when $\rho>0$ is small enough.
\begin{prop}\label{lagrangemultiplier}
Let $2<p<3$ be fixed. For every $\epsilon>0$ there exists
$\rho(\epsilon)>0$ such that
$$\sup_{\omega\in {\mathcal A}_\rho}|\omega - \omega_0|<\epsilon
\hbox{ } \forall 0<\rho<\rho(\epsilon)$$
where
$${\mathcal A}_\rho=
\Big\{\omega\in \R| -\Delta v +\omega v +\rho
\Big (|v|^2*\frac 1{|x|}\Big )v-v|v|^{p-2}=0, v\in {\mathcal N}_{\rho,p}
\Big \}.$$
\end{prop}

{\bf Proof.}
By looking at the equation satisfied by $v\in {\mathcal N}_{\rho,p}$
we deduce
$$\omega=
\frac{\|v\|_{L^p}^p - \|\nabla v\|^2_{L^2}
- \rho \int\int \frac{|v(x)|^2|v(y)|^2}{|x-y|}dxdy}{\|v\|_{L^2}^2}.$$
The proof can be concluded since by Proposition \ref{compactness}
we get that the r.h.s. converges to
$$\frac{\|Q\|_{L^p}^p - \|\nabla Q\|^2_{L^2}}{\|Q\|_{L^2}^2}=\omega_0$$
for $\rho\rightarrow 0$.

\hfill$\Box$

\section{The implicit function argument}

In this section we present some
results strictly related to the
implicit function theorem (see \cite{AP}).
\begin{prop}\label{pro}
There exist $\epsilon_0, \epsilon_1>0$ such that
$$\forall u\in B_{H^1}(Q, \epsilon_0)
\hbox{ } \exists ! \tau(u)\in \R^3,  R(u)\in T^{\perp}
_{Q_{\tau(u)}} \hbox{ s.t. }$$$$
\max\{\|\tau(u)\|_{\R^3}, \|R(u)\|_{H^1}\}<\epsilon_1
\hbox{ and }
u=Q_{\tau(u)}+R(u).
$$
Moreover
$\lim_{u\rightarrow Q} \|\tau(u)\|_{\R^3}=\lim_{u\rightarrow 0} \|R(u)\|_{H^1}=0$
and the nonlinear operators
$$P:B_{H^1}(Q, \epsilon_0)\rightarrow {\mathcal G}$$
$$R:B_{H^1}(Q, \epsilon_0)\rightarrow H^1$$
(where $P(u)=Q_{\tau(u)}$ and $R(u)$ is defined as above) are smooth.
\end{prop}

\begin{remark}\label{radialiP}
Notice that every radially symmetric function
$u\in H^1_{rad}$ can be written as
$u=Q+ (u-Q)$
and moreover $u-Q\in {T_Q^\perp}$ (this follows by noticing
that $T_Q=span \{\partial_{x_i}Q|i=1,...,n\}$).
In particular $Pu=Q$ for every $u\in H^1_{rad}$.
\end{remark}
\begin{remark}\label{traslation}
Notice that $T_{Q_\tau}=\{v(x+\tau)|v\in T_{Q}
\}$. As a consequence it is easy to prove
$P (u(x+\tau))=P(u)(x+\tau)$ and hence (since $P(Q)=Q$)
$P (Q_\tau)=Q_\tau$.
\end{remark}

{\bf Proof.}
It is sufficient to apply the implicit function theorem to the map
$$\Phi:{\mathcal G}\times H^1\ni
(Q_\tau, h)\rightarrow (Q_\tau+h, (h,v_1(x+\tau)),(h,v_2(x+\tau)),
(h, v_3(x+\tau)))\in H^1\times\R^3$$
where $span \{v_1,v_2,v_3\}=T_Q$ and $(,)$ denotes
the usual $L^2$ scalar product.\\
Next we shall prove that
$$d\Phi_{(Q, 0)}\in {\mathcal L} (T_Q\times H^1, H^1\times \R^3)$$ is invertible.
By explicit computation we get
$$d\Phi_{(Q, 0)}:T_Q\times H^1\ni (w, k)\rightarrow
(w+k, (k, v_1), (k, v_2), (k, v_3))\in H^1\times \R^3$$
and hence:
$$\{d\Phi_{(Q, 0)}(-h,h)|h\in T_Q\}=\{0\}\times \R^3;$$
$$\{d\Phi_{(Q, 0)}(0,h)|h\in T_Q^\perp \}= T_Q^\perp \times \{0\};$$
$$\{d\Phi_{(Q, 0)}(h,0)|h\in T_Q \}= T_Q \times \{0\}.$$
As a consequence we deduce that $d\Phi_{(Q, 0)}$ is surjective.\\
Next we prove that $d\Phi_{(Q,0)}$ is injective.
Assume that $(w, k)\in  T_Q \times H^1$ satisfy
$$d\Phi_{(Q, 0)}(w,k)=(0,0)\in H^1\times \R^3$$
which in turn (by looking at the explicit structure of $d\Phi_{(Q,0)}$
) is equivalent to
$$k\in T_Q^\perp
\hbox{ and }
w=-k,$$
hence $w\in T_Q^\perp$. By combining this fact with the hypothesis
$w\in T_Q$ we get $w=0$ and also
$k=-w=0$.

\hfill$\Box$

\begin{prop}\label{proradial}
There exists $\epsilon_2>0$
such that the equation
$$-\Delta w + \omega w +\rho \Big (|w|^2*\frac 1{|x|}\Big )w-w|w|^{p-2}=0$$
has a solution
$w(\rho, \omega)\in H^1_{rad}$
for every $(\rho, \omega)\in (0, \epsilon_2)\times (\omega_0 -\epsilon_2,
\omega_0+\epsilon_2)$.
Moreover
$$\lim_{(\omega, \rho)\rightarrow
(\omega_0, 0)} w(\omega, \rho)= Q \hbox{ in } H^1.
$$
\end{prop}

{\bf Proof.} It follows by an application of
the implicit function theorem at the following operator:
$$\Phi: \R\times \R^+
\times H^1_{rad} \ni (\omega, \rho, u)
\rightarrow \nabla_u {\mathcal F}_{\rho, \omega,p}\in H^1_{rad}$$
where
$$ {\mathcal F}_{\rho,\omega,p}(u)=
{\mathcal E}_{\rho,p}(u)+\frac \omega 2 \|u\|_{L^2}^2.$$
Notice that $\nabla_Q {\mathcal F}_{0,\omega_0, p}=0$
and moreover
\begin{equation}\label{struc}
 d\Phi_{(0, \omega_0, Q)}[h]= h+ K h \hbox{ } \forall h\in H^1_{rad}
\end{equation}
where
$$(p-1)^{-1}K= (-\Delta +\omega_0)^{-1} \circ |Q|^{p-2}.$$
Due to the decay properties of the function $Q(x)$
and the Rellich Compactness Theorem (see \cite{B})
the operator
$H^1_{rad}\ni v\rightarrow |Q|^{p-2} v\in L^2_{rad}$ is compact. Moreover
$(-\Delta +\omega_0)^{-1}\in {\mathcal L}(L^2_{rad}, H^1_{rad})$
and hence the operator $K\in {\mathcal L}(H^1_{rad}, H^1_{rad})$ is a compact operator.
By combining this fact with \eqref{struc} we deduce that
$d\Phi_{(0, \omega_0, Q)}\in {\mathcal L}(H^1_{rad}, H^1_{rad})$
is a Fredholm operator with index zero (see \cite{B}).\\
Moreover by the work \cite{W} it is easy to deduce that
$$\ker_{H^1_{rad}} d\Phi_{(0, \omega_0, Q)} =\{h\in H^1_{rad}(\R^3)| h+ Kh=0\}=\{0\}$$
and hence $d\Phi_{(0, \omega_0, Q)}$ is invertible (since
$d\Phi_{(0, \omega_0, Q)}$ is injective and has Fredholm index zero).

\hfill$\Box$

In next proposition (and along its proof) the operators $P(u)$,$R(u)$
and the number $\epsilon_0>0$ are the ones
in Proposition \ref{pro}.
\begin{prop}\label{proimpo}
There exist $\epsilon_3, \epsilon_4>0$
such that:
\begin{equation}\label{unique}\forall (\omega, \rho)\in (\omega_0-
\epsilon_3,
\omega_0+\epsilon_3)\times (0, \epsilon_3)
\exists ! \hbox{ } u=u(\rho,  \omega)\in H^1
\hbox{ s.t. }\end{equation}
$$\|u(\rho, \omega)-Q\|_{H^1}<\epsilon_4,
P(u)=Q \hbox{ and }
\pi_{T_Q^\perp} (\nabla_u {\mathcal F}_{\rho, \omega,p})=0$$
where
\begin{equation}\label{Fgrande} {\mathcal F}_{\rho,\omega,p}(u)=
{\mathcal E}_{\rho,p}(u)+\frac \omega 2 \|u\|_{L^2}^2.\end{equation}
\end{prop}

{\bf Proof.}
It is sufficient to apply the implicit function theorem to the map
$$\Phi: \R\times \R\times B_{H^1}
(Q, \epsilon_0)\ni (\rho, \omega, u)\rightarrow
(\pi_{T_Q^\perp} (\nabla_u {\mathcal F}_{\rho,\omega,p}),
P(u))\in T_Q^\perp \times {\mathcal G}.
$$
Hence we have to show that
$d\Phi_{(0,\omega_0,Q)}\in {\mathcal L}(H^1,T_{Q}^\perp\times T_Q)$
is invertible.
Recall that by remark \ref{traslation} we get
$P (Q_\tau)=Q_\tau$ and hence
$\Phi(0,\omega_0, Q_\tau)=(0,  Q_\tau)$ which in turn implies
\begin{equation}\label{triv}
d\Phi_{(0,\omega_0,Q)}[v]= (0,v) \hbox{ } \forall v\in T_Q.
\end{equation}
Arguing as in Proposition \ref{proradial}
we deduce that the operator
\begin{equation}\label{fredholm}
d(\nabla_u {\mathcal F}_{\rho, \omega, p})_{(0,\omega_0, Q)}\in {\mathcal L}
(H^1, H^1)\end{equation}$$
\hbox{ is a Fredholm operator of index zero in $H^1$}.$$
Moreover by the work \cite{W} we get
\begin{equation}\label{nondegen}
T_Q=\ker d (\nabla_u {\mathcal F})_{(0, \omega_0, Q)}
\end{equation}
and by the self-adjointness (w.r.t.
to the $L^2$ scalar product) of
the operator $d(\nabla_u {\mathcal F}_{\rho,\omega,p})_{(0,\omega_0, Q)}$
we get
\begin{equation}\label{self}
d(\nabla_u {\mathcal F}_{\rho,\omega,p})_{(0,\omega_0,Q)} (T_Q^\perp)\subset T_Q^\perp.
\end{equation}
By combining \eqref{fredholm}, \eqref{nondegen}
and \eqref{self} we conclude that
\begin{equation}\label{orthiso}d\Phi_{(0,\omega_0,Q)}\in {\mathcal L} (T_Q^\perp, T_Q^\perp)
\hbox{ is invertible .}\end{equation}
By combining \eqref{triv} and
\eqref{orthiso} it is easy to deduce that
$d\Phi_{(0,\omega_0,Q)}\in {\mathcal L}(H^1,T_{Q}^\perp\times T_Q)$
is invertible.

\hfill$\Box$

\section{Proof of proposition \ref{equiv}}

Recall that the operators $P(u),R(u)$ are the ones introduced along
Proposition \ref{pro}.\\
Let
$v\in {\mathcal N}_{\rho,p}$. Due to Proposition \ref{compactness}
for every $\epsilon>0$ there exists $\rho_1(\epsilon)>0$
such that (up-to translation) $v\in B_{H^1}(Q, \epsilon)$ provided that
$\rho<\rho_1(\epsilon)$.
Moreover
$v$ solves the problem
$$-\Delta v + \omega v+\rho \Big
(|v|^2*\frac 1{|x|}\Big ) v- v|v|^{p-2}=0$$
or equivalently
\begin{equation}\label{fullgrad}
\nabla_v {\mathcal F}_{\rho, \omega,p}=0
\end{equation}
(see \eqref{Fgrande} for definition of
${\mathcal F}_{\rho, \omega,p}$) for a suitable $\omega$ such that
$|\omega-\omega_0|<\epsilon$ provided that $\rho<\rho_2(\epsilon)$
(see Proposition \ref{lagrangemultiplier}).
Notice that by Proposition \ref{pro}
there exists $\delta(\epsilon)>0$
such that we can write in a unique way (provided that $\epsilon>0$ is small enough)
$v(x)=Q_\tau + r(x)$ with $\|r\|_{H^1}<\delta$ and $r\in T_{Q_\tau}^\perp$,
and hence
$v(x-\tau)=Q +r(x-\tau)$
with
$r(x-\tau)\in T_Q^\perp$.
By combining this fact with \eqref{fullgrad} (recall the
translation invariance of the functional
${\mathcal F}_{\rho, \omega,p}$) we get
\begin{equation}\label{proper}
P(v(x-\tau))=Q \hbox{ and }
\pi_{T_Q^\perp} (\nabla_{v(x-\tau)} {\mathcal F}_{\rho, \omega,p})=0.
\end{equation}
On the other hand by combining remark \ref{radialiP} with Proposition \ref{proradial}
we deduce that
$w(\rho, \omega)\in H^1_{rad}$ (given in Proposition
\ref{proradial}) satisfies the same properties of $v(x-\tau)$ in \eqref{proper}.
By the uniqueness property included in Proposition \ref{proimpo}
(see \eqref{unique})
we get $v(x-\tau)=w(\rho, \omega)$ and hence $v(x-\tau)$ is radially symmetric.

\end{document}